\documentclass[11pt,reqno]{amsart}
\usepackage{graphicx}
\usepackage{amsfonts}
\usepackage{mathrsfs}
 \usepackage{times}
 \usepackage{mathptmx}
\usepackage{amsmath}
\usepackage{txfonts}
\usepackage{amssymb}
\usepackage{amsmath,amsthm}
\usepackage{amsfonts}
\usepackage[body={7.5in,8.75in}, top=1.2in, right=0.9in,left=0.9in]{geometry}
\usepackage[hypertex,colorlinks=true]{hyperref}
\allowdisplaybreaks[4] \numberwithin{equation}{section}

\newcommand{\abs}[1]{\left\lvert#1\right\rvert}
\newtheorem{thm}{Theorem}

\newtheorem{cor}[thm]{Corollary}
\newtheorem{lem}[thm]{Lemma}
\newtheorem{rem}[thm]{Remark}
\newtheorem{pro}[thm]{Proposition}

\DeclareMathOperator{\td}{d\mspace{-0.5mu}}
\def\squarebox#1{\hbox to #1{\hfill\vbox to #1{\vfill}}}

\begin{document}

\title[On generalized $q$-difference equations  for general Al-Salam--Carlitz polynomials]{A note on generalized $q$-difference equations  for general Al-Salam--Carlitz polynomials}
\author{Jian Cao$^1$, Binbin Xu$^2$ and Sama Arjika${}^{3,*}$ }
\dedicatory{\textsc}
\thanks{${}^{3,*}$Address for correspondence: Department of Mathematics and Informatics,   University of Agadez, Niger. ${}^1$Department
of Mathematics, Hangzhou Normal University, Hangzhou City, Zhejiang Province, 311121, China. ${}^2$Department
of Mathematics, Hangzhou Normal University, Hangzhou City, Zhejiang Province, 311121, China.}
\thanks{Email: 21caojian@hznu.edu.cn,418272337@qq.com, rjksama2008@gmail.com.}



\begin{abstract}
In this paper, we  deduce the generalized $q$-difference equations  for general Al-Salam--Carlitz polynomials and generalize  Arjika's recently results [$q$-difference equation for homogeneous $q$-difference operators and their applications, J. Differ. Equ.  Appl.  {\bf 26}, 987--999 (2020)]. In addition, we obtain transformational identities by the method of $q$-difference equation. Moreover, we deduce $U(n+1)$ type generating functions and Ramanujan's integrals involving general Al-Salam--Carlitz polynomials by $q$-difference equation.
\end{abstract}

\maketitle

\section{Introduction}

In this paper, we refer to the general references \cite{GasparRahman} for definitions and notations. Throughout this paper, we suppose that $0<q<1$. For complex numbers $a$, the $q$-shifted factorials are defined by:
$$ (a;q)_n=\left\{
\begin{aligned}
1,\quad\quad\quad\quad\quad\quad\quad \quad n&=0  \\
(1-a)(1-aq)...(1-aq^{n-1}),\quad n&=1,2,3...
\end{aligned}
\right.
$$
and $(a_1,a_2,...,a_m;q)_n=(a_1;q)_n(a_2;q)_n...(a_m;q)_n$, where $m$ is a positive integer and $n$ is a non-negative integer or $\infty$.\par
    The $q$-binomial coefficient is defined by
\begin{equation}
    \begin{bmatrix}
                        n \\
                       k \\
    \end{bmatrix}=\frac{(q;q)_n}{(q;q)_k(q;q)_{n-k}}.
\end{equation}\par
 
The basic (or $q$-) hypergeometric function 
of the variable $z$ and with $\mathfrak{r}$ numerator 
and $\mathfrak{s}$ denominator parameters 
is defined as follows (see, for details, the monographs by 
\mbox{Slater~(\cite{Slater}, Chapter 3)}
and by Srivastava and Karlsson~(\cite{SrivastavaKarlsson}, 
p. 347, Eq. (272)); 
see also \cite{Koekock,HMS-IMAJAM1983-1984a,
HMS-IMAJAM1983-1984b}): 
 
$${}_{\mathfrak r}\Phi_{\mathfrak s}\left[
\begin{array}{rr}
a_1, a_2,\cdots, a_{\mathfrak r};\\
\\
b_1,b_2,\cdots,b_{\mathfrak s};
\end{array}\,
q;z\right]
:=\sum_{n=0}^\infty\Big[(-1)^n \; 
q^{\binom{n}{2}}\Big]^{1+{\mathfrak s}-{\mathfrak r }}
\,\frac{(a_1, a_2,\cdots, a_{\mathfrak r};q)_n}
{(b_1,b_2,\cdots,b_{\mathfrak s};q)_n}
\; \frac{z^n}{(q;q)_n},$$
where $q\neq 0$ when ${\mathfrak r }>{\mathfrak s}+1$. 
We also note that
$${}_{\mathfrak r+1}\Phi_{\mathfrak r}\left[
\begin{array}{rr}
a_1, a_2,\cdots, a_{\mathfrak r+1}\\
\\
b_1,b_2,\cdots,b_{\mathfrak r };
\end{array}\,
q;z\right]
=\sum_{n=0}^\infty \frac{(a_1, a_2,\cdots, a_{\mathfrak r+1};q)_n}
{(b_1,b_2,\cdots,b_{\mathfrak r};q)_n}\;\frac{ z^n}{(q;q)_n}.$$\par
We remark in passing that, in a recently-published 
survey-cum-expository review article, the~so-called $(p,q)$-calculus 
was exposed to be a rather trivial and inconsequential variation of 
the classical $q$-calculus, the additional parameter $p$ being redundant 
or superfluous \mbox{(see, for details,~(\cite{HMS-ISTT2020}, p. 340)).} 
\par
 
The basic (or $q$-) 
series and basic (or $q$-) polynomials, especially
the basic (or $q$-) hypergeometric functions and basic 
(or $q$-) hypergeometric polynomials, are 
applicable particularly in several diverse areas 
[see also~(\cite{SrivastavaKarlsson}, pp. 350--351)]. In particular, the
celebrated Chu-Vandermonde summation theorem and its known $q$-extensions,
which   have already been 
demonstrated to be useful (see, for details, 
\cite{GasparRahman,Srivastava2012,GEAndrews1986,HMS-JC-SA2020}).

    The usual $q$-differential operator, or $q$-derivative, is defined by \cite{14,15,SrivastaAbdlhusein}
\begin{align}
    D_a\{f(a)\}=\frac{f(a)-f(aq)}{a},\quad
    \theta_a\{f(a)\}=\frac{f(aq^{-1})-f(a)}{q^{-1}a}.
\end{align}\par
    The Leibniz rule for $D_a$ and $\theta_a$ are the following identities \cite{14,15,11}
\begin{align}
    D_a^n\{f(a)g(a)\}&=\sum_{k=0}^nq^{k(k-n)}\begin{bmatrix}n \\k \\ \end{bmatrix}D_a^k\{f(a)\}D_a^{n-k}\Bigl\{g\bigl(aq^k\bigr)\Bigr\},\\
     \theta_a^n\{f(a)g(a)\}&=\sum_{k=0}^n\begin{bmatrix}n \\k \\ \end{bmatrix}\theta_a^k\{f(a)\}\theta_a^{n-k}\Bigl\{g\bigl(aq^{-k}\bigr)\Bigr\}.
\end{align} \par
    The following property of $D_q$ is straightforward and important \cite{18}
\begin{align*}
    D_q\left\{\frac{1}{(at;q)_\infty}\right\}=\frac{t}{(at;q)_\infty},\quad
    D_q^k\left\{\frac{1}{(at;q)_\infty}\right\}=\frac{t^k}{(at;q)_\infty},
\end{align*}
\begin{equation}
D_q^k\{a^n\}=\left\{
\begin{aligned}
\frac{(q,q)_n}{(q,q)_{n-k}}a^{n-k},&\quad0\leq k\leq n-1  \\
(q,q)_n,\quad\quad\quad &\quad k=n \\
0,\quad\quad\quad\quad\quad &\quad k\geq n+1.
\end{aligned}
\right.
\end{equation}\par
    The Al-Salam--Carlitz polynomials were introduced by Al-Salam and Carlitz in 1965 \cite[Eqs. (1.11) and (1.15)]{1}
\begin{align}
   \phi_n^{(a)}(x|q)=\sum_{k=0}^n\begin{bmatrix}
                        n \\
                       k \\
    \end{bmatrix}(a;q)_kx^k\quad and\quad
    \psi_n^{(a)}(x|q)=\sum_{k=0}^n\begin{bmatrix}
                        n \\
                       k \\
    \end{bmatrix}q^{k(k-n)}(aq^{1-k};q)_kx^k.
\end{align}\par
    They play important roles in the theory of $q$-orthogonal polynomials. In fact, there are two families of these polynomials: one with continuous orthogonality and another with discrete orthogonality, which given explicitly in the book of Koekoek--Swarttouw--Lesky {\cite[Eqs. (14.24.1) and (14.25.1)]{5}}. For further information about $q$-polynomials, see \cite{cao-niu2016,cao-sri-liu,cao,5,wang2012,wang2013}.\par
    The generalized Al-Salam--Carlitz polynomials \cite[Eq. (4.7)]{17}
\begin{align}
\phi_n^{(a,b,c)}(x,y|q)&=\sum_{k=0}^n\begin{bmatrix}
                        n \\
                       k \\
    \end{bmatrix}\frac{(a,b;q)_k}{(c;q)_k}x^ky^{n-k},\nonumber\\
\psi_n^{(a,b,c)}(x,y|q)&=\sum_{k=0}^n\begin{bmatrix}
                        n \\
                       k \\
    \end{bmatrix}\frac{(-1)^kq^{{k+1\choose2}-nk}(a,b;q)_k}{(c;q)_k}x^ky^{n-k},
\end{align}
 whose generating functions are \cite[Eqs. (4.10) and (4.11))]{17}
\begin{align}
\sum_{n=0}^\infty \phi_n^{(a,b,c)}(x,y|q)\frac{t^n}{(q;q)_n}&=\frac{1}{(yt;q)_\infty}{}_2\phi_1\Biggl[\begin{matrix}
\begin{array}{r}
a,b;\\\\
c;
\end{array}
\end{matrix}q;xt\Biggr],\quad \max\{\abs{yt},\abs{xt}\}<1,\label{9}\\
\sum_{n=0}^\infty \psi_n^{(a,b,c)}(x,y|q)\frac{(-1)^nq^{n\choose 2}t^n}{(q;q)_n}&=(yt;q)_\infty{}_2\phi_1\Biggl[\begin{matrix}
\begin{array}{r}
a,b;\\\\
c;
\end{array}
\end{matrix}q;xt\Biggr],\quad \abs{xt}<1.\label{10}
\end{align}
\par
Chen and Liu \cite{14,15} gave the clever way of parameter augmentation by use of the following two $q$-exponential operators
\begin{align}
   \mathbb{T}(bD_a)=\sum_{n=0}^\infty\frac{(bD_a)^n}{(q;q)_n},\quad
   \mathbb{E}(b\theta_a)=\sum_{n=0}^\infty\frac{q^{n\choose2}(b\theta_a)^n}{(q;q)_n},
\end{align}
which is rich and powerful tool for basic hypergeometric series, especially makes many famous results easily fall into this framework. For further information about $q$-exponential operators, see \cite{14,15,fang2013,fang2014,fang2015,jia2014}.\par
Recently, Srivastava, Arjika and Sherif Kelil \cite{6}  introduced the following homogeneous $q$-difference operator $\widetilde{E}(a,b;D_q)$ by
\begin{equation}\label{Srivastava-Arjika-Kelil }
   \widetilde{E}(a,b;D_q)=\sum_{n=0}^\infty\frac{(-1)^nq^{n\choose2}(a;q)_n}{(q;q)_n}(bD_a)^n.
\end{equation}
The operators \eqref{Srivastava-Arjika-Kelil } have turned out to be suitable for dealing with a generalized Cauchy polynomials $p_n(x,y,a)$
\begin{equation}
  p_n(x,y,a)=\widetilde{E}(a,b;D_q)\{x^n\}.
\end{equation}
For more information about the relations between operators and $q$-polynomials, see \cite{6}.\par
Liu \cite{18,16} deduced several results involving Bailey's $_6\psi_6$, $q$-Mehler formulas for Rogers--Szeg\"{o} polynomials and $q$-integral of Sears' transformation by the following $q$-difference equations.

\begin{pro}[{\cite[Theorems 1 and 2]{16}}]\label{pro1}
Let $f(a,b)$ be a two-variable analytic function in a neighbourhood of $(a,b)=(0,0)\in\mathbb{C}^2$\par
(I) If $f(a,b)$ satisfied the difference equation
\begin{equation}
bf(aq,b)-af(a,bq)=(b-a)f(a,b),
\end{equation}\par
\quad\quad then we have
\begin{equation}
f(a,b)=\mathbb{T}(bD_a)\{f(a,0)\}.
\end{equation}\par
(II) If $f(a,b)$ satisfied the difference equation
\begin{equation}
af(aq,b)-bf(a,bq)=(a-b)f(aq,bq),
\end{equation}\par
\quad\quad then we have
\begin{equation}
f(a,b)=\mathbb{E}(b\theta_a)\{f(a,0)\}.
\end{equation}
\end{pro}
 Arjika \cite{10} continue to consider the following generalized $q$-difference equations.
\begin{pro}[{\cite[Theorem 2.1]{10}}]\label{pro3}
Let $f(a,x,y)$ be a three-variable analytic function in a neighbourhood of $(a,x,y)=(0,0,0)\in\mathbb{C}^3$. If $f(a,x,y)$ can be expanded in terms of $p_n(x,y,a)$ if and only if
\begin{equation}\label{pro3_1}
\quad\quad\quad x[f(a,x,y)-f(a,x,qy)]=y[f(a,qx,qy)-f(a,x,qy)]-ay[f(a,qx,q^2y)-f(a,x,q^2y)].
\end{equation}
\end{pro}
\begin{pro}[{\cite[Theorem 2.2]{10}}]\label{pro2}
Let $f(a,x,y)$ be a three-variable analytic function in a neighbourhood of $(a,x,y)=(0,0,0)\in\mathbb{C}^3$.
If $f(a,x,y)$ satisfied the $q$-difference equation
\begin{equation}\label{pro2_1}
\quad\quad\quad x[f(a,x,y)-f(a,x,qy)]=y[f(a,qx,qy)-f(a,x,qy)]-ay[f(a,qx,q^2y)-f(a,x,q^2y)],
\end{equation}\par
then we have
\begin{equation}
f(a,x,y)=\widetilde{E}(a,b;D_q)\{f(a,x,0)\}.
\end{equation}
\end{pro}
In this paper, our goal is to generalize the results of   Arjika \cite{10}  in section 2. We first construct the following $q$-operators
\begin{align}\label{1.9}
   \mathbb{T}(a,b,c,d,e,yD_x)&=\sum_{n=0}^\infty\frac{(a,b,c;q)_n}{(q,d,e;q)_n}(yD_x)^n,\\
   \mathbb{E}(a,b,c,d,e,y\theta_x)&=\sum_{n=0}^\infty\frac{(-1)^nq^{n\choose2}(a,b,c;q)_n}{(q,d,e;q)_n}(y\theta_x)^n.
\end{align}

We remark that the $q$-operator (\ref{1.9}) is a particular case of  
the homogeneous $q$-difference operator $\mathbb{T}({\bf a},{\bf b},cD_{x})$  
(see \cite{HMS-Sama2020}) by taking  
$${\bf a} =(a,b,c),\quad {\bf b} =(d,e) \qquad \text{and} \qquad c=y.$$

We also built the relations between operators $\mathbb{T}(a,b,c,d,e,yD_x)$, $\mathbb{E}(a,b,c,d,e,y\theta_x)$ and  the new generalized Al-Salam--Carlitz polynomials $\phi_n^{a,b,c\choose d,e}(x,y|q)$, $\psi_n^{a,b,c\choose d,e}(x,y|q),$ respectively,
\begin{align}
\phi_n^{a,b,c\choose d,e}(x,y|q)\triangleq\mathbb{T}(a,b,c,d,e,yD_x)\{x^n\}&=\sum_{k=0}^n\begin{bmatrix}
                        n \\
                       k \\
    \end{bmatrix}\frac{(a,b,c;q)_k}{(d,e;q)_k}x^{n-k}y^k,\label{23}\\
\psi_n^{a,b,c\choose d,e}(x,y|q)\triangleq\mathbb{E}(a,b,c,d,e,y\theta_x)\{x^n\}&=\sum_{k=0}^n\begin{bmatrix}
                        n \\
                       k \\
    \end{bmatrix}\frac{(-1)^kq^{k(k-n)}(a,b,c;q)_k}{(d,e;q)_k}x^{n-k}y^k.\label{24}
\end{align} 

The paper is organized as follows: In section \ref{section2}, we state two theorems and give the proofs. In section \ref{section3}, we gain generalize generating functions for new generalized Al-Salam--Carlitz polynomials by using the method of $q$-difference equations perspectively. In section \ref{section4}, we obtain a transformational identities involving generating functions for generalized Al-Salam--Carlitz polynomials by $q$-difference equations. In section \ref{section5}, we deduce $U(n+1)$ type generating functions for generalized Al-Salam--Carlitz polynomials by $q$-difference equation. In section \ref{section6}, we deduce generalizations of Ramanujan's integrals.

\section{Main results and proofs}\label{section2}

In this section, we give the following two theorems.
\begin{thm}\label{thm1}
Let $f(a,b,c,d,e,x,y)$ be a seven-variable analytic function in a neighbourhood of $(a,b,c,d,e,x,y)=(0,0,0,0,0,0,0)\in\mathbb{C}^7$.\par
(I) If $f(a,b,c,d,e,x,y)$ can be expanded in terms of $\phi_n^{a,b,c\choose d,e}(x,y|q)$ if and only if
\begin{align}\label{thm1_1}
&\quad\quad x\{f(a,b,c,d,e,x,y)-f(a,b,c,d,e,x,yq)-(d+e)q^{-1}[f(a,b,c,d,e,x,yq)-f(a,b,c,d,e,x,yq^2)]\nonumber\\
&\quad\quad\quad +deq^{-2}[f(a,b,c,d,e,x,yq^2)-f(a,b,c,d,e,x,yq^3)]\}\nonumber\\
&=y\{[f(a,b,c,d,e,x,y)-f(a,b,c,d,e,xq,y)]-(a+b+c)[f(a,b,c,d,e,x,yq)-f(a,b,c,d,e,xq,yq)]\nonumber\\
&\quad\quad\quad +(ab+ac+bc)[f(a,b,c,d,e,x,yq^2)-f(a,b,c,d,e,xq,yq^2)]\nonumber\\
&\quad\quad\quad\quad  -abc[f(a,b,c,d,e,x,yq^3)-f(a,b,c,d,e,xq,yq^3)]\}.
\end{align}\par
(II) If $f(a,b,c,d,e,x,y)$ can be expanded in terms of $\psi_n^{a,b,c\choose d,e}(x,y|q)$ if and only if
\begin{align}\label{thm1_2}
&\quad\quad x\{f(a,b,c,d,e,xq,y)-f(a,b,c,d,e,xq,yq)-(d+e)q^{-1}[f(a,b,c,d,e,xq,yq)-f(a,b,c,d,e,xq,yq^2)]\nonumber\\
&\quad\quad\quad +deq^{-2}[f(a,b,c,d,e,xq,yq^2)-f(a,b,c,d,e,xq,yq^3)]\}\nonumber\\
&=y\{[f(a,b,c,d,e,xq,yq)-f(a,b,c,d,e,x,yq)]-(a+b+c)[f(a,b,c,d,e,xq,yq^2)-f(a,b,c,d,e,x,yq^2)]\nonumber\\
&\quad\quad\quad +(ab+ac+bc)[f(a,b,c,d,e,xq,yq^3)-f(a,b,c,d,e,x,yq^3)]\nonumber\\
&\quad\quad\quad\quad  -abc[f(a,b,c,d,e,xq,yq^4)-f(a,b,c,d,e,x,yq^4)]\}.
\end{align}
\end{thm}

\begin{rem}
For $c$=$d$=$e$=$0$, and $b$$\rightarrow$ $\frac{1}{b}$,  $y\rightarrow yb$, $b\rightarrow 0$, then equation \eqref{thm1_1}  reduces  \eqref{pro3_1}.
\end{rem}

\begin{thm}\label{thm2}
Let $f(a,b,c,d,e,x,y)$ be a seven-variable analytic function in a neighbourhood of $(a,b,c,d,e,x,y)=(0,0,0,0,0,0,0)\in\mathbb{C}^7$.

(I) If $f(a,b,c,d,e,x,y)$ satisfied the difference equation
\begin{align}\label{thm2_1}
&\quad\quad x\{f(a,b,c,d,e,x,y)-f(a,b,c,d,e,x,yq)-(d+e)q^{-1}[f(a,b,c,d,e,x,yq)-f(a,b,c,d,e,x,yq^2)]\nonumber\\
&\quad\quad\quad +deq^{-2}[f(a,b,c,d,e,x,yq^2)-f(a,b,c,d,e,x,yq^3)]\}\nonumber\\
&=y\{[f(a,b,c,d,e,x,y)-f(a,b,c,d,e,xq,y)]-(a+b+c)[f(a,b,c,d,e,x,yq)-f(a,b,c,d,e,xq,yq)]\nonumber\\
&\quad\quad\quad +(ab+ac+bc)[f(a,b,c,d,e,x,yq^2)-f(a,b,c,d,e,xq,yq^2)]\nonumber\\
&\quad\quad\quad\quad  -abc[f(a,b,c,d,e,x,yq^3)-f(a,b,c,d,e,xq,yq^3)]\},
\end{align}
then we have
\begin{align}\label{thm2_1.2}
f(a,b,c,d,e,x,y)=\mathbb{T}(a,b,c,d,e,yD_x)\{f(a,b,c,d,e,x,0)\}.
\end{align}
(II)  If $f(a,b,c,d,e,x,y)$ satisfied the difference equation
\begin{align}\label{thm2_2}
&\quad\quad x\{f(a,b,c,d,e,xq,y)-f(a,b,c,d,e,xq,yq)-(d+e)q^{-1}[f(a,b,c,d,e,xq,yq)-f(a,b,c,d,e,xq,yq^2)]\nonumber\\
&\quad\quad\quad +deq^{-2}[f(a,b,c,d,e,xq,yq^2)-f(a,b,c,d,e,xq,yq^3)]\}\nonumber\\
&=y\{[f(a,b,c,d,e,xq,yq)-f(a,b,c,d,e,x,yq)]-(a+b+c)[f(a,b,c,d,e,xq,yq^2)-f(a,b,c,d,e,x,yq^2)]\nonumber\\
&\quad\quad\quad +(ab+ac+bc)[f(a,b,c,d,e,xq,yq^3)-f(a,b,c,d,e,x,yq^3)]\nonumber\\
&\quad\quad\quad\quad  -abc[f(a,b,c,d,e,xq,yq^4)-f(a,b,c,d,e,x,yq^4)]\},
\end{align}
then we have
\begin{align}
f(a,b,c,d,e,x,y)=\mathbb{E}(a,b,c,d,e,y\theta_x)\{f(a,b,c,d,e,x,0)\}.
\end{align}
\end{thm}

\begin{rem}
For $c$=$d$=$e$=$0$, and $b$$\rightarrow$ $\frac{1}{b}$,  $y\rightarrow yb$, $b\rightarrow 0$, then equation \eqref{thm2_1} reduces to \eqref{pro2_1}.
\end{rem}

To determine if a given function is an analytic function in several complex variables, we often use the following Hartogs's theorem. For more information, please refer to \cite{19,7}.

\begin{lem}[{{\cite[Hartogs's theorem]{9}}}]
If a complex-valued function is holomorphic (analytic) in each variable separately in an open domain $D \in\mathbb{C}^n$, then it is holomorphic (analytic) in $D$.
\end{lem}

In order to prove Theorem \ref{thm1}, we need the following fundamental property of several complex variables.
\begin{lem}[{\cite[Proposition 1]{3}}]
If $f(x_1,x_2,...,x_k)$ is analytic at the origin $(0,0,...,0)\in\mathbb{C}^k$, then, $f$ can be expanded in an absolutely convergent power series,
\begin{align*}
f(x_1,x_2,...,x_k)=\sum_{n_1,n_2,...,n_k=0}^\infty\alpha_{n_1,n_2,...,n_k}x_1^{n_1}x_2^{n_2}...x_k^{n_k}.
\end{align*}
\end{lem}

\begin{proof}[Proof of Theorem \ref{thm1}]
(I) From the Hartogs's theorem and the theory of several complex variables, we assume that
\begin{align}\label{27}
f(a,b,c,d,e,x,y)=\sum_{k=0}^\infty A_k(a,b,c,d,e,x)y^k.
\end{align}\par
On one hand, substituting \eqref{27} into \eqref{thm1_1} yields
\begin{multline}
 x\sum_{k=0}^\infty[1-q^k-(d+e)q^{k-1}+(d+e)q^{2k-1}+deq^{2k-2}-deq^{3k-2}]A_k(a,b,c,d,e,x)y^k\\
=\sum_{k=0}^\infty[1-(a+b+c)q^k+(ab+bc+ac)q^{2k}-abcq^{3k}][A_k(a,b,c,d,e,x)-A_k(a,b,c,d,e,xq)]y^{k+1},
\end{multline}
which is equal to
\begin{multline}\label{A-n rel-eq}
 x\sum_{k=0}^\infty(1-q^k)(1-dq^{k-1})(1-eq^{k-1})A_k(a,b,c,d,e,x)y^k\\
=\sum_{k=0}^\infty(1-aq^{k})(1-bq^{k})(1-cq^{k})[A_k(a,b,c,d,e,x)-A_k(a,b,c,d,e,xq)]y^{k+1}.
\end{multline}\par
Equating coefficients of $y^k$ on both sides of equation \eqref{A-n rel-eq}, we have
\begin{multline}
 x(1-q^k)(1-dq^{k-1})(1-eq^{k-1})A_k(a,b,c,d,e,x)\\
=(1-aq^{k-1})(1-bq^{k-1})(1-cq^{k-1})[A_{k-1}(a,b,c,d,e,x)-A_{k-1}(a,b,c,d,e,xq)],
\end{multline}
which is equivalent to
\begin{align*}
A_k(a,b,c,d,e,x)&=\frac{(1-aq^{k-1})(1-bq^{k-1})(1-cq^{k-1})}{(1-q^k)(1-dq^{k-1})(1-eq^{k-1})}\cdot\frac{A_{k-1}(a,b,c,d,e,x)-A_{k-1}(a,b,c,d,e,xq)}{x} \\
&=\frac{(1-aq^{k-1})(1-bq^{k-1})(1-cq^{k-1})}{(1-q^k)(1-dq^{k-1})(1-eq^{k-1})}\cdot D_x\{A_{k-1}(a,b,c,d,e,x)\}.
\end{align*}\par
By iteration, we gain
\begin{align*}
A_k(a,b,c,d,e,x)=\frac{(a,b,c;q)_k}{(q,d,e;q)_k}\cdot D_x^k\{A_0(a,b,c,d,e,x)\}.
\end{align*}\par
Letting $\displaystyle f(a,b,c,d,e,x,0)=A_0(a,b,c,d,e,x)=\sum_{n=0}^\infty\mu_nx^n$ yields
\begin{align}\label{28}
A_k(a,b,c,d,e,x)=\frac{(a,b,c;q)_k}{(q,d,e;q)_k}\cdot\sum_{n=0}^\infty\mu_n\frac{(q;q)_n}{(q;q)_{n-k}}x^{n-k},
\end{align}
we have
\begin{align*}
f(a,b,c,d,e,x,y)&=\sum_{k=0}^\infty\frac{(a,b,c;q)_k}{(q,d,e;q)_k}\sum_{n=0}^\infty\mu_n\frac{(q;q)_n}{(q;q)_{n-k}}x^{n-k}y^k\\
&=\sum_{n=0}^\infty\mu_n\sum_{k=0}^\infty\begin{bmatrix}
                        n \\
                       k \\
    \end{bmatrix}\frac{(a,b,c;q)_k}{(d,e;q)_k}x^{n-k}y^k\\
&=\sum_{n=0}^\infty\mu_n\phi_n^{a,b,c\choose d,e}(x,y|q).
\end{align*}
On the other hand, if $f(a,b,c,d,e,x,y)$ can be expanded in terms of $\phi_n^{a,b,c\choose d,e}(x,y|q)$, we verify that $f(a,b,c,d,e,x,y)$ satisfies \eqref{thm1_1}. Similarly, we deduce (II). The proof of Theorem \ref{thm1} is complete.
\end{proof}

\begin{proof}[Proof of Theorem \ref{thm2}]
From the theory of several complex variables, we begin to solve the $q$-difference. First we may assume that
\begin{align}\label{33}
f(a,b,c,d,e,x,y)=\sum_{k=0}^\infty A_k(a,b,c,d,e,x)y^k.
\end{align}\par
Substituting this equation into \eqref{33} and compare coefficients of $y^k\,\,(k\geq1)$, we readily find that
\begin{multline}
 x(1-q^k)(1-dq^{k-1})(1-eq^{k-1})A_k(a,b,c,d,e,x)\\
=(1-aq^{k-1})(1-bq^{k-1})(1-cq^{k-1})[A_{k-1}(a,b,c,d,e,x)-A_{k-1}(a,b,c,d,e,xq)],
\end{multline}
which equals
\begin{align*}
A_k(a,b,c,d,e,x)&=\frac{(1-aq^{k-1})(1-bq^{k-1})(1-cq^{k-1})}{(1-q^k)(1-dq^{k-1})(1-eq^{k-1})}\cdot\frac{A_{k-1}(a,b,c,d,e,x)-A_{k-1}(a,b,c,d,e,xq)}{x} \\
&=\frac{(1-aq^{k-1})(1-bq^{k-1})(1-cq^{k-1})}{(1-q^k)(1-dq^{k-1})(1-eq^{k-1})}\cdot D_x\{A_{k-1}(a,b,c,d,e,x)\}.
\end{align*}\par
By iteration, we gain
\begin{align}\label{34}
A_k(a,b,c,d,e,x)=\frac{(a,b,c;q)_k}{(q,d,e;q)_k}\cdot D_x^k\{A_0(a,b,c,d,e,x)\}.
\end{align}\par
Now we return to calculate $A_0(a,b,c,d,e,x)$. Just taking $y=0$ in equation \eqref{33}, we immediately obtain $A_0(a,b,c,d,e,x)=f(a,b,c,d,e,x,0)$. Substituting \eqref{34} into \eqref{33}, we achieve \eqref{thm2_1.2} directly. The proof of Theorem \ref{thm2} is complete.
\end{proof}


\section{Generating functions for new generalized Al-Salam--Carlitz polynomials}\label{section3}
In this section we generalized generating functions for the new generalized Al-Salam--Carlitz polynomials by the method of $q$-difference equations.\par
We start with the following lemmas.
\begin{lem}[\cite{13}]
The Cauchy polynomials as given below
\begin{align}\label{lem_3}
p_n(x,y)=(x-y)(x-qy)...(x-q^{n-1}y)=(y/x;q)_nx^n
\end{align}
together with the following Srivastava-Agarwal 
type generating function 
(see also \cite{Cao-Srivastava2013}):  
\begin{equation}
\label{Srivas}
\sum_{n=0}^\infty  
p_n (x,y)\;\frac{(\lambda;q)_n\,t^n}{(q;q)_n}
={}_{2}\Phi_1\left[
\begin{array}{rr}
\lambda,\frac{y}{x};\\
\\
0;
\end{array}\,
q;  xt\right].
\end{equation}
\end{lem}
\begin{lem}[\cite{13}]\label{lem_4}
Suppose that $max\{|xt|,|yt|\}<1$, we have
\begin{align}
\sum_{n=0}^\infty p_n(x,y)\frac{t^n}{(q;q)_n}=\frac{(yt;q)_\infty}{(xt;q)_\infty}.\label{putt}
\end{align}
\end{lem}

Based upon the $q$-binomial theorem or   the homogeneous version of the Cauchy identity
(\ref{putt}) and Heine's
transformations, Srivastava {et al.} \cite{HMS-C-W2020} established 
a set of two presumably new theta-function identities 
(see, for details, \cite{HMS-C-W2020}).
\begin{lem}[{\cite[Theorem 5]{13}}]
Suppose that $\max\{\abs{act},\abs{adt},\abs{bct},\abs{bdt}\}<1$, we have
\begin{align}\label{lem_5}
\sum_{n=0}^\infty h_n(a,b|q)h_n(c,d|q)\frac{t^n}{(q;q)_n}=\frac{(abcdt^2;q)_\infty}{(act,adt,bct,bdt;q)_\infty}.
\end{align}
\end{lem}
\begin{thm}\label{thm3}
Suppose that $max\{|xt|,|yt|\}<1$, we have
\begin{align}\label{thm3_1}
\sum_{n=0}^\infty \phi_n^{a,b,c\choose d,e}(x,y|q)\frac{t^n}{(q;q)_n}=\frac{1}{(xt;q)_\infty}{}_3\phi_2\Biggl[\begin{matrix}
\begin{array}{r}
a,b,c;\\\\
d,e;
\end{array}
\end{matrix}q;yt\Biggr],
\end{align}
\begin{align}\label{thm3_2}
\sum_{n=0}^\infty \psi_n^{a,b,c\choose d,e}(x,y|q)\frac{t^n}{(q;q)_n}=(xt;q)_\infty{}_3\phi_3\Biggl[\begin{matrix}
\begin{array}{r}
a,b,c;\\\\
0,d,e;
\end{array}
\end{matrix}q;-yt\Biggr].
\end{align}
\end{thm}
\begin{rem}
For $c=e=0$ in Theorem \ref{thm3}, equations \eqref{thm3_1} and \eqref{thm3_2} reduce to equations \eqref{9} and \eqref{10}, respectively.
\end{rem}
\begin{proof}[Proof of Theorem \ref{thm3}]
By denoting the right-hand side of equation \eqref{thm3_1} by $f(a,b,c,d,e,x,y)$,we can verify that $f(a,b,c,d,e,x,y)$ satisfies \eqref{thm1_1}. So, we have
\begin{align*}
f(a,b,c,d,e,x,y)=\sum_{k=0}^\infty\mu_n\phi_n^{a,b,c\choose d,e}(x,y|q),
\end{align*}
and
\begin{align*}
f(a,b,c,d,e,x,0)=\sum_{k=0}^\infty\mu_nx^n=\frac{1}{(xt;q)_\infty}=\sum_{n=0}^\infty\frac{t^n}{(q;q)_n}x^n.
\end{align*}
So, $f(a,b,c,d,e,x,y)$ is equal to
\begin{align*}
f(a,b,c,d,e,x,y)=\sum_{k=0}^\infty\frac{t^n}{(q;q)_n}\phi_n^{a,b,c\choose d,e}(x,y|q),
\end{align*}
equal to the right-hand side of equation \eqref{thm3_1}. \par
Similarly, by denoting the right-hand side of equation \eqref{thm3_2} by $f(a,b,c,d,e,x,y)$, we can verify that $f(a,b,c,d,e,x,y)$ satisfies \eqref{thm1_2}. So, we can use the same way to achieve the equation \eqref{thm3_2}. The proof of Theorem \ref{thm3} is complete.
\end{proof}

\begin{thm}\label{thm4}
Suppose that $\max\{\abs{xt},\abs{yt}\}<1$, we have
\begin{align}\label{thm4_1}
\sum_{n=0}^\infty \phi_n^{a,b,c\choose d,e}(x,y|q)\frac{p_n(s,t)}{(q;q)_n}=\frac{(xs;q)_\infty}{(xt;q)_\infty}{}_4\phi_3\Biggl[\begin{matrix}
\begin{array}{r}
a,b,c,s/t;\\\\
d,e,xs;
\end{array}
\end{matrix}q;yt\Biggr].
\end{align}
\end{thm}
\begin{cor}
Suppose that $\abs{yt}<1$, we have
\begin{align}\label{cor1_1}
\sum_{n=0}^\infty \phi_n^{a,b,c\choose d,e}(x,y|q)\frac{(-1)^nq^{n\choose2}t^n}{(q;q)_n}=(xt;q)_\infty{}_3\phi_3\Biggl[\begin{matrix}
\begin{array}{r}
a,b,c;\\\\
d,e,xt;
\end{array}
\end{matrix}q;yt\Biggr].
\end{align}
\end{cor}
\begin{rem}
For $t=0$, in Theorem \ref{thm4}, equation \eqref{thm4_1} reduces to \eqref{cor1_1}. For $s=0$ in Theorem \ref{thm4}, equation \eqref{thm4_1} reduces to \eqref{thm3_1}, respectively.
\end{rem}

\begin{proof}[Proof of Theorem \ref{thm4}]
By denoting the right-hand side of equation \eqref{thm4_1} by $f(a,b,c,d,e,x,y)$, we can verify that $f(a,b,c,d,e,x,y)$ satisfies \eqref{thm1_1}. So, we have
\begin{align}
f(a,b,c,d,e,x,y)=\sum_{k=0}^\infty\mu_n\phi_n^{a,b,c\choose d,e}(x,y|q),
\end{align}
and
\begin{align*}
f(a,b,c,d,e,x,0)=\sum_{k=0}^\infty\mu_nx^n=\frac{(xs;q)_\infty}{(xt;q)_\infty}=\sum_{n=0}^\infty\frac{p_n(t,s)}{(q;q)_n}x^n.
\end{align*}
So, $f(a,b,c,d,e,x,y)$ is equal to the right-hand side of equation \eqref{thm4_1}. The proof of Theorem \ref{thm4} is complete.
\end{proof}

\begin{thm}\label{thm5}
For $k\in\mathbb{N}$ and $max\{|xt|,|yt|\}<1$, we have
\begin{align}\label{thm5_1}
\sum_{n=0}^\infty \phi_{n+k}^{a,b,c\choose d,e}(x,y|q)\frac{t^n}{(q;q)_n}=\frac{x^k}{(xt;q)_\infty}\sum_{n=0}^\infty\frac{(a,b,c;q)_n}{(q,d,e;q)_n}(yt)^n\sum_{j=0}^n\begin{bmatrix}
                        n \\
                       k \\
    \end{bmatrix}\frac{(-1)^jq^{nj-{j\choose2}}(q^{-k},xt;q)_j}{(xt)^j}.
\end{align}
\end{thm}
\begin{rem}
For $k=0$ in Theorem \ref{thm5}, equation \eqref{thm5_1} reduces to \eqref{thm3_1}.
\end{rem}
\begin{proof}[Proof of Theorem \ref{thm5}]
Denoting the right-hand side of equation \eqref{thm5_1} equivalently by
\begin{align}\label{44}
f(a,b,c,d,e,x,y)=\frac{x^k}{(xt;q)_\infty}\sum_{n=0}^\infty\frac{(a,b,c;q)_n}{(q,d,e;q)_n}(yt)^n\sum_{j=0}^n\begin{bmatrix}
                        n \\
                       k \\
    \end{bmatrix}\frac{(-1)^jq^{nj-{j\choose2}}(q^{-k},xt;q)_j}{(xt)^j},
\end{align}
and it is easy to check that \eqref{44} satisfies \eqref{thm1_1}, so we have
\begin{align}\label{45}
f(a,b,c,d,e,x,y)=\sum_{k=0}^\infty\mu_n\phi_n^{a,b,c\choose d,e}(x,y|q).
\end{align}
Setting $y=0$ in \eqref{45}, it becomes
\begin{align*}
f(a,b,c,d,e,x,0)=\sum_{k=0}^\infty\mu_nx^n=\frac{x^n}{(xt;q)_\infty}=\sum_{n=0}^\infty x^k\frac{(xt)^n}{(q;q)_n}=\sum_{n=0}^\infty x^{n+k}\frac{t^n}{(q;q)_n}=\sum_{n=k}^\infty x^n\frac{t^{n-k}}{(q;q)_{n-k}}.
\end{align*}
Hence
\begin{align*}
f(a,b,c,d,e,x,y)=\sum_{n=k}^\infty\phi_n^{a,b,c\choose d,e}(x,y|q)\frac{t^{n-k}}{(q;q)_{n-k}}=\sum_{n=0}^\infty \phi_{n+k}^{a,b,c\choose d,e}(x,y|q)\frac{t^n}{(q;q)_n}.
\end{align*}\par
The proof of Theorem \ref{thm5} is complete.
\end{proof}
\begin{thm}\label{thm6}
We have
\begin{multline}\label{thm6_1}
 \sum_{n=0}^\infty \phi_n^{a_1,b_1,c_1\choose d_1,e_1}(x_1,y_1|q)\phi_n^{a_2,b_2,c_2\choose d_2,e_2}(x_2,y_2|q)\frac{t^n}{(q;q)_n}=\frac{1}{(x_1x_2t;q)_\infty}\sum_{n=0}^\infty\frac{(a_2,b_2,c_2;q)_n(x_1y_2t)^n}{(q,d_2,e_2;q)_n}\\
\times\sum_{j=0}^\infty\frac{(q^{n-j+1},x_1x_2t,a_1,b_1,c_1;q)_j(\frac{y_1}{x_1})^j}{(q,d_1,e_1;q)_j}
{}_3\phi_2\Biggl[\begin{matrix}
\begin{array}{r}
a_1q^j,b_1q^j,c_1q^j;\\\\
d_1q^j,e_1q^j;
\end{array}
\end{matrix}q;x_2y_1t\Biggr].
\end{multline}
\end{thm}
\begin{rem}
For $a_1$=$b_1$=$c_1$=$d_1$=$e_1$=$a_2$=$b_2$=$c_2$=$d_2$=$e_2$=$0$ in Theorem \ref{thm6}, equation \eqref{thm6_1} reduces to \eqref{lem_5}.
\end{rem}
\begin{proof}[Proof of Theorem \ref{thm6}]
Denoting the right-hand side of equation \eqref{thm6_1} by $H(a_1,b_1,c_1,d_1,e_1,x_1,y_1)$, we have
\begin{multline}\label{47}
\quad H(a_1,b_1,c_1,d_1,e_1,x_1,y_1)=\frac{1}{(x_1x_2t;q)_\infty}\sum_{n=0}^\infty\frac{(a_2,b_2,c_2;q)_n(x_1y_2t)^n}{(q,d_2,e_2;q)_n}\\
\times\sum_{j=0}^\infty\frac{(q^{n-j+1},x_1x_2t,a_1,b_1,c_1;q)_j(\frac{y_1}{x_1})^j}{(q,d_1,e_1;q)_j}
{}_3\phi_2\Biggl[\begin{matrix}
\begin{array}{r}
a_1q^j,b_1q^j,c_1q^j;\\\\
d_1q^j,e_1q^j;
\end{array}
\end{matrix}q;x_2y_1t\Biggr].
\end{multline}\par
Because equation \eqref{47} satisfies \eqref{thm2_1}, we have
\begin{align*}
H(a_1,b_1,c_1,d_1,e_1,x_1,y_1)&=\mathbb{T}(a_1,b_1,c_1,d_1,e_1,y_1D_{x_1})\{H(a_1,b_1,c_1,d_1,e_1,x_1,0)\}\nonumber\\
&=\mathbb{T}(a_1,b_1,c_1,d_1,e_1,y_1D_{x_1})\left\{\frac{1}{(x_1x_2t;q)_\infty}\sum_{n=0}^\infty\frac{(a_2,b_2,c_2;q)_n(x_1y_2t)^n}{(q,d_2,e_2;q)_n}\right\}\\
&=\mathbb{T}(a_1,b_1,c_1,d_1,e_1,y_1D_{x_1})\left\{\sum_{n=0}^\infty\phi_n^{a_2,b_2,c_2\choose d_2,e_2}(x_2,y_2|q)\frac{(x_1t)^n}{(q;q)_n}\right\}\\
&=\sum_{n=0}^\infty\phi_n^{a_2,b_2,c_2\choose d_2,e_2}(x_2,y_2|q)\frac{t^n}{(q;q)_n}\mathbb{T}(a_1,b_1,c_1,d_1,e_1,y_1D_{x_1})\{x_1^n\}\\
&=\sum_{n=0}^\infty\phi_n^{a_1,b_1,c_1\choose d_1,e_1}(x_1,y_1|q) \phi_n^{a_2,b_2,c_2\choose d_2,e_2}(x_2,y_2|q)\frac{t^n}{(q;q)_n}.
\end{align*}\par
The proof of Theorem \ref{thm6} is complete.
\end{proof}

\section{Transformational identities from $q$-difference equations }\label{section4}
Liu \cite{17} gave some important transformational identities by the method of $q$-difference operator. For more details, please refer to \cite{5,17,20}.\par
In this section we deduce the following transformational identities involving generating functions for new generalized Al-Salam--Carlitz polynomials by the method of $q$-difference equation.
\begin{thm}\label{thm7}
Let $A(k)$ and $B(k)$ satisfy
\begin{align}\label{thm7_1}
\sum_{k=0}^\infty A(k)x^k=\sum_{k=0}^\infty B(k)\frac{(xtq^k;q)_\infty}{(xq^k;q)_\infty},
\end{align}
and we have
\begin{align}\label{thm7_2}
\sum_{k=0}^\infty A(k)\phi_k^{a,b,c\choose d,e}(x,y|q)=\sum_{k=0}^\infty B(k)\frac{(xtq^k;q)_\infty}{(xq^k;q)_\infty}{}_4\phi_3\Biggl[\begin{matrix}
\begin{array}{r}
a,b,c,1/t;\\\\
d,e,xtq^k;
\end{array}
\end{matrix}q;yq^k\Biggr],
\end{align}
supposing that equations \eqref{thm7_1} and \eqref{thm7_2} are convergent.
\end{thm}
\begin{cor}\label{cor2}
Suppose that $|r|,|x|,|xt|<1$, we have
\begin{align}\label{cor2_1}
\sum_{k=0}^\infty\phi_k^{a,b,c\choose d,e}(x,y|q)\frac{(t,s;q)_k}{(q,r;q)_k}=\frac{(xt,s;q)_\infty}{(x,r;q)_\infty}\sum_{k=0}^\infty\frac{(r/s,x;q)_ks^k}{(q,xt;q)_k}{}_4\phi_3\Biggl[\begin{matrix}
\begin{array}{r}
a,b,c,1/t;\\\\
d,e,xtq^k;
\end{array}
\end{matrix}q;yq^k\Biggr].
\end{align}
\end{cor}
\begin{rem}
Setting $A(k)$ and $B(k)$ in Theorem \ref{thm7} by \eqref{53} given below, equation \eqref{thm7_2} reduces to \eqref{cor2_1}. For $y=0$ in \eqref{cor2_1}, equation \eqref{thm7_2}  reduces to \eqref{52} below.
\end{rem}
\begin{proof}[Proof of Theorem \ref{thm7}]
Denoting the right-hand side of equation \eqref{thm7_2} equivalently by $f(a,b,c,d,e,x,y)$, and we can check that $f(a,b,c,d,e,x,y)$ satisfied \eqref{thm1_1}, so we have
\begin{align}\label{51}
f(a,b,c,d,e,x,y)=\sum_{k=0}^\infty\mu_n\phi_n^{a,b,c\choose d,e}(x,y|q)
\end{align}\\
Setting $y=0$ in \eqref{51}, it becomes
\begin{align*}
f(a,b,c,d,e,x,0)&=\sum_{k=0}^\infty\mu_nx^n=\sum_{k=0}^\infty B(k)\frac{(xtq^k;q)_\infty}{(xq^k;q)_\infty}\quad by\,\,(48)\\
&=\sum_{k=0}^\infty A(k)x^k.
\end{align*}
Hence
\begin{align*}
f(a,b,c,d,e,x,y)=\sum_{k=0}^\infty A(k)\phi_k^{a,b,c\choose d,e}(x,y|q).
\end{align*}\par
The proof of Theorem \ref{thm7} is complete.
\end{proof}
\begin{proof}[Proof of Corollary \ref{cor2}]
Using the Heine's $q$-Euler transformations \cite[Eq.(1.4.1)]{1}
\begin{align}\label{52}
{}_2\phi_1\Biggl[\begin{matrix}\begin{array}{r}t,s;\\\\r;\end{array}\end{matrix}q;x\Biggr]=\frac{(s,xt;q)_\infty}{(r,x;q)_\infty}
{}_2\phi_1\Biggl[\begin{matrix}\begin{array}{r}r/s,x;\\\\xt;\end{array}\end{matrix}q;s\Biggr],
\end{align}
formula \eqref{thm7_1} is valid for the case
\begin{align}\label{53}
A(k)=\sum_{k=0}^\infty\frac{(t,s;q)_k}{(q,r)_k}\quad and \quad B(k)=\frac{(s;q)_\infty}{(r;q)_\infty}\sum_{k=0}^\infty\frac{(r/s;q)_k}{(q;q)_k}s^k.
\end{align}\par
Using equation \eqref{thm7_2}, we can deduce the Corollary \ref{cor2}.
\end{proof}

\section{$U(n+1)$ type generating functions for generalized Al-Salam-Carlitz polynomials }\label{section5}
    Multiple basic hypergeometric series associated to the unitary $U(n+1)$  group have been investigated by various authors, see \cite{8,12}. In \cite{8}, Milne initiated theory and application of the $U(n+1)$ generalization of the classical Bailey transform and Bailey lemma, which involves the following nonterminating $U(n+1)$  generalizations of the $q$-binomial theorem.
\begin{pro}[{\cite[Theorem 5.42]{18}}]
Let $b,z$ and $x_1,...,x_n$ be indeterminate, and let $n\geq 1$. Suppose that none of the denominators in the following identity vanishes, and that $0<|q|<1$ and $|z|<|x_1,...,x_n||x_m|^{-n}|q|^{(n-1)/2}$, for $m=1,2,...,n$. Then
\begin{multline}\label{lem6_1}
\sum_{\mbox{\tiny$\begin{array}{c} y_n \geq 0\\ k=1,2,...,n\end{array}$}}\biggl\{\prod_{1\leq r<s\leq n}\biggl[\frac{1-(x_r/x_sq^{y_r-y_s})}{1-(x_r/x_s)}\biggl]\prod_{r,s=1}^n(q\frac{x_r}{x_s};q)^{-1}_{y_r}\prod_{i=1}^n(x_i)^{ny_i-(y_1+...+y_n)}(-1)^{(n-1)(y_1+...+y_n)}\\
\times q^{y_2+2y_3+...+(n-1)y_n+(n-1)[{y_1\choose 2}+...+{y_n\choose 2}]-e_2(y_1,...,y_n)}(b;q)_{y_1+...+y_n}z^{y_1+...+y_n}=\frac{(bz;q)_\infty}{(z;q)_\infty},
\end{multline}
where $e_2(y_1,...,y_n)$ is the second elementary symmetric function of $\{y_1,...,y_n\}$.
\end{pro}

In this section, we deduce $U(n+1)$ type generating functions for generalized Al-Salam-Carlitz polynomials by the methodn of $q$-difference equation.
\begin{thm}\label{thm8}
Let $b,z$ and $x_1,...,x_n$ be indeterminate, and let $n\geq 1$. Suppose that none of the denominators in the following identity vanishes, and that $0<|q|<1$ and $|z|<|x_1,...,x_n||x_m|^{-n}|q|^{(n-1)/2}$, for $m=1,2,...,n$. Then
\begin{multline}\label{thm8_1}
\sum_{\mbox{\tiny$\begin{array}{c} y_n \geq 0\\ k=1,2,...,n\end{array}$}}\biggl\{\prod_{1\leq r<s\leq n}\biggl[\frac{1-(x_r/x_sq^{y_r-y_s})}{1-(x_r/x_s)}\biggl]\prod_{r,s=1}^n(q\frac{x_r}{x_s};q)^{-1}_{y_r}\prod_{i=1}^n(x_i)^{ny_i-(y_1+...+y_n)}(-1)^{(n-1)(y_1+...+y_n)}\\
\times q^{y_2+2y_3+...+(n-1)y_n+(n-1)[{y_1\choose 2}+...+{y_n\choose 2}]-e_2(y_1,...,y_n)}\phi_{y_1+...+y_n}^{r,s,t\choose u,v}(z,y|q)(b;q)_{y_1+...+y_n}=\frac{(bz;q)_\infty}{(z;q)_\infty}{}_4\phi_3\Biggl[\begin{matrix}
\begin{array}{r}
r,s,t,b;\\\\
u,v,bz;
\end{array}
\end{matrix}q;y\Biggr],
\end{multline}
where $e_2(y_1,...,y_n)$ is the second elementary symmetric function of $\{y_1,...,y_n\}$.
\end{thm}
\begin{rem}
For $y=0$, in Theorem \ref{thm8}, equation \eqref{thm8_1} reduces to \eqref{lem6_1}.
\end{rem}
\begin{proof}[Proof of Theorem \ref{thm8}]
Denoting the right-hand side of equation \eqref{thm8_1} equivalently by $f(r,s,t,u,v,z,y)$, and we can check that $f(r,s,t,u,v,z,y)$ satisfied \eqref{thm2_1}), so we have
\begin{align*}
f(r,s,t,u,v,z,y)&=\mathbb{T}(r,s,t,u,v,yD_z)\{f(r,s,t,u,v,z,0)\}\\
&=\sum_{\mbox{\tiny$\begin{array}{c} y_n \geq 0\\ k=1,2,...,n\end{array}$}}\biggl\{\prod_{1\leq r<s\leq n}\biggl[\frac{1-(x_r/x_sq^{y_r-y_s})}{1-(x_r/x_s)}\biggl]\prod_{r,s=1}^n(q\frac{x_r}{x_s};q)^{-1}_{y_r}\prod_{i=1}^n(x_i)^{ny_i-(y_1+...+y_n)}(-1)^{(n-1)(y_1+...+y_n)}\\
&\quad\times q^{y_2+2y_3+...+(n-1)y_n+(n-1)[{y_1\choose 2}+...+{y_n\choose 2}]-e_2(y_1,...,y_n)}(b;q)_{y_1+...+y_n}\mathbb{T}(r,s,t,u,v,yD_z)\{z^{y_1+...+y_n}\},
\end{align*}
which is the left-hand side of \eqref{thm8_1} by \eqref{23}. The proof of Theorem \ref{thm8} is complete.
\end{proof}

\section{Generalization of Ramanujan's integrals }\label{section6}
The following integral of Ramanujan \cite{2} are quite famous.
\begin{pro}[{\cite[Eq. (2)]{2}}]\label{lem7}
For $0<q=e^{-2k^2}<1$ and $m\in\mathbb{R}$. Suppose that $\abs{abq}<1$, we have
\begin{align}\label{lem7_1}
\int^\infty_\infty\frac{e^{-x^2+2mx}}{(aq^{1/2}e^{2ikx},bq^{1/2}e^{-2ikx};q)_\infty}\td x=\sqrt{\pi}e^{m^2}\frac{(-aqe^{2mki},-bqe^{-2mki};q)_\infty}{(abq;q)_\infty}.
\end{align}
\end{pro}
In this section, we have the following generalization of Ramanujan's integrals.
\begin{thm}\label{thm9}
For $0<q=e^{-2k^2}<1$ and $m\in\mathbb{R}$. Suppose that $\abs{abq}<1$, we have
\begin{multline}\label{thm9_1}
\int^\infty_\infty\frac{e^{-x^2+2mx}}{(aq^{1/2}e^{2ikx},bq^{1/2}e^{-2ikx};q)_\infty}{}_3\phi_2\Biggl[\begin{matrix}
\begin{array}{r}
r,s,t;\\\\
u,v;
\end{array}
\end{matrix}q;yq^{1/2}e^{2ikx}\Biggr]\td x\\
=\sqrt{\pi}e^{m^2}\frac{(-aqe^{2mki},-bqe^{-2mki};q)_\infty}{(abq;q)_\infty}{}_4\phi_3\Biggl[\begin{matrix}
\begin{array}{r}
r,s,t,e^{2mki}/b;\\\\
u,v,-aqe^{2mki};
\end{array}
\end{matrix}q;ybq\Biggr].
\end{multline}
\end{thm}
\begin{rem}
For $y=0$ in Theorem \ref{thm9}, equation \eqref{thm9_1} reduces to \eqref{lem7_1}.
\end{rem}
\begin{proof}[Proof of Theorem \ref{thm9}]
Denoting the right-hand side of \eqref{thm9_1} equivalently by $f(r,s,t,u,v,a,y)$. $f(r,s,t,u,v,a,y)$ is analytic near $(r,s,t,u,v,a,y)$ and we can check that $f(r,s,t,u,v,a,y)$ satisfied \eqref{thm1_1}, so we have
\begin{align*}
f(r,s,t,u,v,a,y)=\sum_{k=0}^\infty\mu_n\phi_n^{r,s,t\choose u,v}(a,y|q),
\end{align*}\\
and
\begin{align*}
f(r,s,t,u,v,a,0)&=\sum_{k=0}^\infty\mu_na^n=\sqrt{\pi}e^{m^2}\frac{(-aqe^{2mki},-bqe^{-2mki};q)_\infty}{(abq;q)_\infty}\quad\, by\,\,\,\eqref{lem7_1}\\
&=\int^\infty_\infty\frac{e^{-x^2+2mx}}{(aq^{1/2}e^{2ikx},bq^{1/2}e^{-2ikx};q)_\infty}\td x\\
&=\int^\infty_\infty\frac{e^{-x^2+2mx}}{(bq^{1/2}e^{-2ikx};q)_\infty}\left\{\sum_{n=0}^\infty\frac{(q^{1/2}e^{2ikx})^n}{(q;q)_n}a^n\right\}\td x.
\end{align*}\par
So we have
\begin{align*}
f(r,s,t,u,v,a,y)=\int^\infty_\infty\frac{e^{-x^2+2mx}}{(bq^{1/2}e^{-2ikx};q)_\infty}\left\{\sum_{n=0}^\infty\phi_n^{r,s,t\choose u,v}(a,y|q)\frac{(q^{1/2}e^{2ikx})^n}{(q;q)_n}\right\}\td x,
\end{align*}\\
which is equal to the left-hand side of equation \eqref{thm9_1}. The proof of Theorem \ref{thm9} is complete.
\end{proof}

\section{ Concluding Remarks and Observations}
\label{conclusion}
In our present investigation, we have introduced a set of 
two $q$-operators \break $\mathbb{T}(a,b,c,d,e,yD_x)$ and 
$\mathbb{E}(a,b,c,d,e,y\theta_x)$ with to applying them
to  generalize Arjika's recently results \cite{10}, and derive transformational identities by means of the $q$-difference equations. We have   also  derived   
$U(n+1)$-type generating functions and Ramanujan's integrals involving general Al-Salam-Carlitz polynomials by means of the $q$-difference equations.

It is believed that the $q$-series and $q$-integral identities, which we
have presented in this paper, as~well as the various related recent works
cited here, will provide encouragement and motivation for further researches
on the topics that are dealt with and investigated in this paper.

In conclusion, we find it to be worthwhile to remark that some potential
further applications of the methodology and findings, which we have
presented here by means of the $q$-analysis and the $q$-calculus, 
can be found in the study of
the zeta and $q$-zeta functions as well as their related 
functions of Analytic Number Theory
(see, for example, \cite{HMS-JAEC2019, HMS-PIMM2019}; see also
\cite{Srivastava2012}) and also in the study of analytic and
univalent functions of Geometric Function Theory {\it via} 
number-theoretic entities (see, for example, \cite{HMS-Symmstry2020}).

\medskip

\noindent
{\bf Conflicts of Interest:} The authors declare that they 
have no conflicts of interest.\\
{\bf Availability of data and material:} Not applicable. \\
 {\bf Funding:} This work was supported by the Zhejiang Provincial 
Natural Science Foundation of China (No. LY21A010019). \\
 {\bf Authors' contributions:} Both authors equally contributed to this manuscript   and approved the final version.  \\
 {\bf Competing interests:}  The authors declare that there are no competing interests.  \\
 {\bf Acknowledgements:} Not applicable.

\end{document}